\def \version {2013--9--23}

\documentclass[12pt]{article}
\usepackage{amssymb}

\def \sxs {simplexes}
\def \sst {\subset}

\def \err {\mathbb{R}}

\def \ex {\mathrm{ex}}
\def \ybl {YBLM}

\def \kkk {{\LARGE\sf \ \mbox{---!!!---} \ }}

\newtheorem{Theorem}{Theorem}
\newtheorem{lem}[Theorem]{Lemma}
\newtheorem{crl}[Theorem]{Corollary}
\newtheorem{prp}[Theorem]{Proposition}
\newtheorem{dfn}[Theorem]{Definition}
\newtheorem{rmk}[Theorem]{Remark}
\newtheorem{xmp}[Theorem]{Example}
\newtheorem{clm}[Theorem]{Claim}
\newtheorem{op}[Theorem]{Problem}

\newtheorem{problem}[Theorem]{Problem}
\newtheorem{conjecture}[Theorem]{Conjecture}

\def \bp {\begin{prp} \ }
\def \ep {\end{prp}}

\def \bd {\begin{dfn} \ }
\def \ed {\end{dfn}}

\def \bc {\begin{crl} \ }
\def \ec {\end{crl}}

\def \bcj {\begin{conjecture} \ }
\def \ecj {\end{conjecture}}

\def \thm {\begin{Theorem} \ }
\def \ethm {\end{Theorem}}

\def \bl {\begin{lem} \ }
\def \el {\end{lem}}

\def \brm {\begin{rmk} \ }
\def \erm {\end{rmk}}

\def \bxm {\begin{xmp} \ \rm }
\def \exm {\end{xmp}}

\def \bcm {\begin{clm} \ }
\def \ecm {\end{clm}}

\def \bop {\begin{op} \ }
\def \eop {\end{op}}

\def \tmz {\begin{itemize}}
\def \etmz {\end{itemize}}

\def \nin {\noindent}
\def \bsk {\bigskip}

\def \smin {\setminus}
\def \es {\emptyset}

\def\cH{{\cal H}}

\def\cE{{\cal E}}
\def\cF{{\cal F}}

\def\cK{{\cal K}}
\def\cS{{\cal S}}

\begin{document}

\title{Minimum Number of Affine Simplexes \\
  of Given Dimension}
\author{Istv\'{a}n Szalkai \thanks{~corresponding author}
 $^,$\thanks{~Department of Mathematics,
University of Pannonia, Veszpr\'{e}m, Hungary; \,e-mail:
\texttt{szalkai@almos.uni-pannon.hu}} \and
Zsolt Tuza
 \thanks{~Alfr\'{e}d R\'{e}nyi Institute of Mathematics,
Hungarian Academy of Sciences, Budapest,
and
Department of Computer Science and Systems Technology,
University of Pannonia, Veszpr\'{e}m, Hungary; \,e-mail:
\texttt{tuza@dcs.uni-pannon.hu} }
 $^,$\thanks{~Research supported in part by the Hungarian Scientific Research
  Fund, OTKA grant T-81493, and
  co-financed
 by the European Social Fund through the project T\'AMOP-4.2.2.C-11/1/KONV-2012-0004 -- National Research Center
 for Development and Market Introduction of Advanced Information and Communication
 Technologies.} }
\date{\small Latest update on \version}
\maketitle

\begin{abstract}
In this paper we formulate and solve extremal problems in the Euclidean space
 $\mathbb{R}^{d}$ and further in hypergraphs, originating from problems in
stoichiometry and elementary linear algebra.
 The notion of affine simplex is the bridge between the original
  problems and the presented extremal theorem on set systems.
 As a sample corollary, it follows that if no triple is collinear
  in a set $S$ of $n$ points in $\mathbb{R}^{3}$, then $S$ contains
  at least ${n\choose 4}-cn^3$ affine simplexes for some constant $c$.
 A function related to Sperner's theorem and the YBLM inequality is also
  considered and its relation to hypergraph Tur\'an problems is discussed.

%
%

\bsk


 \nin

\bsk

 \nin

\end{abstract}

\vfil

\section{Introduction}

 The roots of the present study date back to the mid-1980s,
  to the paper by S. Kumar and
 \'{A}.~Peth\H{o} \cite{KP-ICE1985},
  concerning an application of linear algebra in stoichiometry.
From the algebraic point of view, their very natural question asks about
 the number of those subsets of a set of vectors which are linearly
 dependent but each of whose proper subsets is independent.
Here we give an asymptotically tight solution for the minimum
 in terms of dimension and the number of vectors when low-dimensional
 dependencies are excluded.
Our method is to prove a more general result in extremal set theory
 (Theorem~\ref{n-choo-k} below), hence without
 assuming anything about the structure of algebraic dependencies.

\subsection{Motivation in chemistry}

Restricting attention to a ``universe'' of $D$ kinds of atoms
 (or atomic parts), each molecule (species) can be represented with
 a vector in $\mathbb{R}^D$ whose $i$th coordinate means the number
 of atoms of $i$th type in the molecule in question.\footnote{~The
  types of atoms are supposed to be in a fixed order. E.g., if $D=3$
  and the universe is [C,\,H,\,O], then we have the vector $(0,2,1)$
  for H$_2$O and $(2,4,2)$ for CH$_3$COOH.}
Then a chemical reaction naturally corresponds to a zero-sum
 linear combination of these vectors
 (using the law of mass balance).

The reaction is called minimal if none of the molecules,
 taking role in it, can be omitted so that the remaining ones could
 form still a(nother) reaction.
In the language of linear algebra this assumption is equivalent to
 the property that the corresponding set of vectors is linearly dependent
  but each proper subset of it is independent, that is
 the defining condition of \emph{linear algebraic simplex}.
 Both from practical and theoretical purposes the following
problem was raised:

\begin{problem}
\label{min-max-simp(m)_Probl}What is the minimum and maximum number of
linear algebraic simplexes\/ $S\subseteqq \mathcal{V}$
in a set\/ $\mathcal{V}$ of vectors in\/ $\mathbb{R}^{D}$ if only the size\/ $%
\left\vert \mathcal{V}\right\vert $ is given and\/ $\mathcal{V}$ spans\/ $%
\mathbb{R}^{D}$? What are the structures of sets $\mathcal{V}$ which
contain extremal number of simplexes?
\end{problem}

The answer was given in \cite{SzL-HJIC1995}. Moreover, Problem~\ref%
{min-max-simp(m)_Probl} was generalized for matroids in
 \cite{SzLD-PuMA2006}; actually its authors solved it a decade
 earlier than published, see \cite{DLSz-ms}.

Concerning minimum, the results in \cite{SzL-HJIC1995} show that
almost all vectors must be \textit{parallel\/}, i.e.\ almost all molecules
(species) are isomer molecules or multiple doses. The problem where
parallel vectors are excluded is still unsolved in general:

\begin{problem}     \label{min_simp(m)_alt_Probl}
  What is the minimum number of linear algebraic
simplexes\/ $S\subseteqq \mathcal{V}$ if only the size\/ $\left\vert \mathcal{V}%
\right\vert $ is given,\/ $\mathcal{V}$ does not contain parallel
vectors and\/ $\mathcal{V}$ spans\/ $\mathbb{R}^D$? \ What are the structures
of sets\/ $\mathcal{V}$ which contain the minimum number of simplexes?
\end{problem}

A conjecture on both the minimum number and the structure attaining it
 is stated in \cite{SzL-EJC1998}.
The cases $D=3$ and $D=4$ were solved in \cite{SzL-EJC1998} and
 \cite{SzSz-JMC2011}, respectively.



\subsection{Geometric formulation}

In the framework of linear algebra the problem is somewhat
 non-symmetric because the zero vector plays a special role.
This asymmetry can be eliminated if we translate the problem
 to the language of geometry.
Moreover, restricting attention to sets $\mathcal{V}\subset\mathbb{R}^D$
 containing neither the zero vector nor a pair of parallel vectors,
 dimension can be reduced from $D$ to $d=D-1$:
 first associate each element $\underline{v}\in \mathcal{V}$ with
its direction $\lambda \cdot \underline{v}$
 $\left( \lambda \in \mathbb{R}\right)$, and then
 intersect this system $\Lambda \mathcal{V}$ with a $(D-1)$-dimensional
 hyperplane $\mathcal{P}$ which does not contain the origin and
  is not parallel to any element of $\mathcal{V}$.

 The mapping from $\mathcal{V}$ to the set
 $\mathcal{V}^{\mathcal{P}}:=\Lambda \mathcal{V}\cap \mathcal{P}$\ is a
bijection under which linear algebraic simplexes $S\subset \mathcal{V}$
 correspond to \textit{affine simplexes}
  $S^{\mathcal{P}}\subset \mathbb{R}^{D-1}$, where
 a set $S$ of $k\ge 3$ points in the Euclidean $d$-space is called an
 \emph{affine simplex} if $S$ is contained in some
 $(k-2)$-dimensional hyperplane but no proper subset $S'\subsetneqq S$
 is contained in a hyperplane of dimension $|S'|-2$.
For instance, in $\mathbb{R}^3$ the following three types of
 affine simplexes occur:
 \begin{itemize}
  \item three collinear points;
  \item four coplanar points, no three of which are collinear;
  \item five points, no four of which are coplanar.
 \end{itemize}
Affine simplexes can alternatively be defined by requiring
 that the vectors
 $\underline{s}_{2}-\underline{s}_{1},
  \underline{s}_{3}-\underline{s}_{1},\dots,
   \underline{s}_{k}-\underline{s}_{1}$
  be linearly dependent but their proper subsets shouldn't
  (for every choice of a point to be labeled $\underline{s}_{1}$).

In cases of low dimension, as solved in \cite{SzL-EJC1998} and
 \cite{SzSz-JMC2011}, almost all points of the extremal
  configurations for Problem \ref{min_simp(m)_alt_Probl}
   attaining the minimum number of affine simplexes
 lie on one or two lines, i.e.\ mostly contain affine simplexes
 of three points.
In this way the natural question arises to determine the minimum in the
 other extreme, where no three points are collinear.
For this reason our goal is to study point sets which
 contain no affine simplexes smaller than a given size.
The first interesting case is $\err^3$.

Let $S\sst\err^d$ be a set of $n$ points, no $d$ of which lie on a
 $(d-2)$-dimensional hyperplane.
Then two kinds of subsets of $S$ form an affine simplex:
 \tmz
  \item 
   $d+1$ points on a hyperplane of
    dimension $d-1$, or
  \item 
   $d+2$ points, no $d+1$ of which lie
   on a common hyperplane of dimension $d-1$.
 \etmz

\thm \label{alg-d}
 For every\/ $d\ge 3$ there is a constant\/ $c=c(d)$
  with the following property.
 If\/ $S\sst\err^d$ is a set of\/ $n$ points,
  no\/ $d$ of them lying on a hyperplane
  of dimension\/ $d-2$, then\/ $S$ determines
  at least\/ ${n\choose d+1} - c n^d$ affine \sxs.
\ethm

\bc \label{alg-3}
 For any\/ $n$ points in the 3-space, no three being collinear,
  the number of coplanar quadruples plus
  the 5-tuples containing no coplanar quadruples is at least\/
  ${n\choose 4} - O(n^3)$ as\/ $n\to\infty$.
\ec

These results are asymptotically tight, as shown by the obvious
 example of $n$ coplanar points in $\err^3$
 (no three of them being on a line)
  and also for any $d\ge 3$
 by $n$ points of $\err^{d-1}$ in general position when embedded
 isometrically into $\err^d$.
Such a set of points has exactly ${n\choose d+1}$
 affine simplexes.
In fact, configurations with even fewer affine simplexes
 exist, which in addition span the $d$-space.
For instance, $n-1$ points of $\err^{d-1}$ in general position
 embedded in a hyperplane of $\err^d$ plus an $n$th point outside that
 hyperplane generate just ${n-1\choose d+1}$ affine simplexes
 (as no affine simplex contains the $n$th point).

In $\err^3$, the two arrangements of points just mentioned yield
 $\frac{1}{24}n^4-\frac{1}{4}n^3+O(n^2)$ and
  $\frac{1}{24}n^4-\frac{5}{12}n^3+O(n^2)$, respectively.
Currently we do not know whether or not the latter error term
 $\frac{5}{12}n^3$ is asymptotically tight.
We do know, however, that the construction above is not extremal;
 an improvement of the order $O(n^2)$ will be proved
  in Proposition \ref{lower}.

\subsection{Combinatorial formulation}

Here we put the problems and results above in a more general
setting. Let $\cH=(X,\cE)$ be a hypergraph, where $X$ is the
finite vertex set and $\cE$ is the edge set consisting of subsets
of $X$. We extend the notion of linear hypergraph (also called
``simple'' or ``almost disjoint'' in some parts of the literature)
as follows.

\bd
 We say that a hypergraph\/ $\cH=(X,\cE)$ is\/ \emph{$q$-linear}
 (for some integer\/ $q\ge 1$) if\/ $|E\cap E'| < q$ holds
 for all\/ $E,E'\in\cE$, $E\neq E'$.
\ed

Hence, in a 1-linear hypergraph any two edges are disjoint, and
 2-linear coincides with linear hypergraphs in the usual sense,
  in analogy with Euclidean spaces where
   any two points uniquely determine a line.

We also introduce some notation. As usual, $S\choose k$ will
stand for the collection of all $k$-element subsets of set $S$.
For any hypergraph $\cH=(X,\cE)$, let
 \tmz
  \item $\cE_k := {\bigcup_{E\in\cE}} {E\choose k}$ --- usually
   $(X,\cE_k)$ is called the \emph{$k$-section hypergraph} of $\cH$;
  \item $\cE^0_{k+1} := \{ F\in {X\choose k+1} \mid
   {F\choose k} \cap \cE_k = \es \}$.
 \etmz
Corresponding to $k=d+1$, in analogy with the geometric
interpretation, we call the members of $\cE_k \cup \cE^0_{k+1}$
the $(k-1)$-dimensional \emph{semi-\sxs}\ in~$\cH$.

\thm \label{n-choo-k}
 For every\/ $k\ge 3$ there is a constant\/ $c=c(k)$ such that
  $$|\cE_k|+|\cE^0_{k+1}| \ge {n\choose k} - cn^{k-1}$$ holds
  for all\/ $(k-1)$-linear hypergraphs\/ $\cH=(X,\cE)$ on\/
  $n$ vertices.
\ethm

This result implies Theorem \ref{alg-d}, by considering the hypergraph
 whose edges are the sets of points lying on a common hyperplane of
 dimension $d-1$.

\subsection{Sperner families and Tur\'an numbers}
   \label{s:Sp}

For any hypergraph $\cH=(X,\cE)$ (not necessarily $q$-linear for
 a prescribed value of $q$) and for any $k$, the set system
  $\cS_k(\cH):=\cE_k\cup\cE^0_{k+1}$
 is a \emph{Sperner family}, which means that none of its members
 contains any other $S\in\cS$.

The well known \ybl\ inequality\footnote{~For several decades, it was
 called LYM inequality, stated and proved in exactly that form independently by
 Yamamoto \cite{Y}, Meshalkin \cite{M} and Lubell \cite{L} (in this order of
 chronology). Bollob\'as \cite{B} proved a more general result, however,
 from which the inequality follows immediately. Inequalities of
 this kind have lots of applications in extremal problems in various
 areas of mathematics; cf.\ the two-part survey \cite{T1,T2}.
 The current acronym \ybl\ coincides (apart from punctuation) with the
 abbreviated name of famous Hungarian architect Mikl\'os Ybl
 (1814--1891).}
 states that
 \begin{equation}   \label{Sperner}
   \sum_{S\in\cS} {n\choose |S|}^{-1} \le 1
 \end{equation}
  holds for every Sperner family $\cS$
 (where $n$ is the number of vertices).
In particular, (\ref{Sperner}) is valid for the family
 $\cS=\cS_k(\cH)$ of any $\cH$, too.
In connection with the main problem studied here, one may also
 consider the values
 $$
   s(n,k) :=
    \min_{\cH\mbox{\scriptsize\rm ~is~}(k-1)\mbox{\scriptsize\rm -linear},~|X|=n}
      \sum_{S\in\cS_k(\cH)} {n\choose |S|}^{-1}
 $$
  and analogously, without assuming $(k-1)$-linearity,
 $$
   s'(n,k) :=
    \min_{\cH=(X,\cE),~|X|=n}
      \sum_{S\in\cS_k(\cH)} {n\choose |S|}^{-1}.
 $$
Since there exist only finitely many hypergraphs on any given number $n$
 of vertices, both $s(n,k)$ and $s'(n,k)$ are well-defined and
  are at most 1 by the \ybl\ inequality, for all $n$ and $k$.
In Theorem~\ref{Sper-th} we prove that for every fixed $k$, the values of
 $s(n,k)$ and $s'(n,k)$ tend to constants
  larger than 0 and smaller than 1 as $n$ gets large.
We also consider their relation to the Tur\'an problem
 on graphs and uniform hypergraphs.

\section{Proof of the general lower bound}

Here we prove Theorem \ref{n-choo-k}. By the free choice of
$c=c(k)$, we may restrict ourselves to $n$ sufficiently large, say
$n > k^3$. Moreover, due to the nature of the problem, we may
also assume without loss of generality that $|E|\ge k$ holds for
all $E\in \cE$. Let $H\in {X\choose k}$ be any $k$-tuple. If it is
contained in some $E\in\cE$, then $H$ is counted in $\cE_k$
precisely once. We will prove that, with possibly few exceptions,
also the other $k$-tuples $H$ generate at least one member of
$\cE^0_{k+1}$ on the average.
More explicitly, it will turn out that most of those sets $H$
 can be completed to a member of $\cE^0_{k+1}$ in more than $k$
 different ways.

From now on we assume that $H\in {X\choose k}\smin \cE_k$.
 Let $x\in H$ be any vertex.
If the subset $H\smin \{x\}$ is contained in an
edge of $\cH$, we denote that edge by $E_x$; and otherwise we
define $E_x := H\smin \{x\}$.
 A more precise and unambiguous notation would be $E_x(H)$,
  but for simplicity we write $E_x$ as long as just one $H$
  is considered.
 Note that $E_x$ is unique for each $x\in H$ (once $H$ is
 understood), since $\cH$ is $(k-1)$-linear.
It also follows for any two distinct $x,x'\in H$ that
 $E_x$ and $E_{x'}$ share no vertex outside $H$.
 We set
$$
  H^* := \bigcup_{x\in H} E_x .
$$
Then we have the implication
\begin{equation} \label{csillag}
  H \in {X\choose k}\smin \cE_k \ \land \ z \in X \smin H^*
    \quad \Longrightarrow \quad
  H \cup \{z\} \in \cE^0_{k+1}
\end{equation}
because the containment relation $(H \cup \{z\})\smin\{x\}\sst E'$
for some $x\in H$ and $E'\in\cE$ would contradict the assumption
$|E_x\cap E'| < k-1$.

We say that the $k$-tuple $H$ is a \emph{near-cover} of $\cH$ if
$|H^*| \ge |X|-k$. The proof now splits into two situations,
whether $\cH$ has, or does not have, a near-cover.

Suppose first that no $H\in {X\choose k}\smin \cE_k$ is a
near-cover of $\cH$. Then by (\ref{csillag}) we obtain that each
$H$ can be extended to a member $F$ of $\cE^0_{k+1}$ in at least
$k+1$ different ways. On the other hand, each $F\in \cE^0_{k+1}$
can be obtained from exactly $k+1$ sets $H\in {X\choose k}\smin
\cE_k$, namely from its $k$-element subsets. Thus, in this case we
have
$$
  \left| \cE^0_{k+1} \right| \ge
    \left| {X\choose k}\smin \cE_k \right|
$$
and the inequality stated in the theorem
 holds even without the error term $O(n^{k-1})$.

Suppose now that some $H\in {X\choose k}\smin \cE_k$ is a
 near-cover of $\cH$. Then the cardinality of the set
$$
  X' := H^* \smin H
$$
is at least $n-2k$, and $X'$ is partitioned into sets of type
$E'_x := E_x\smin H$ ($x\in H$). Say, $X' = E'_{x_1}\cup\, \cdots \,
\cup E'_{x_\ell}$ where $\ell \le k$.

A case that can directly be settled is when some
  $\ell-1$ sets from $\{ E'_{x_1}, \dots, E'_{x_\ell} \}$
  cover together at most $k^3-4k^2+5k$ vertices.
There are at most $2k$ vertices outside $X'$, hence some
 $E'_{x_i}$ contains at least $n-k^3+4k^2-7k$ elements.
Then
 we  obtain that
  $|\cE_k| \ge {n-k^3+4k^2-7k \choose k} = {n\choose k} - O(n^{k-1})$
 is valid\footnote{~We may actually write the somewhat larger value
  ${n-k^3+4k^2-7k+k-1 \choose k}$,
 by considering $E_{x_i}$ instead of $E'_{x_i}$; but this is
 irrelevant concerning the current proof.}.
Therefore, we may suppose for the rest of the proof that the union
 of any $\ell-1$ sets from $\{ E'_{x_1}, \dots, E'_{x_\ell} \}$
 has cardinality greater than $k^3-4k^2+5k$.

 Consider any $H'\in {X'\choose k}\smin (\cE_k \cup \{H\})$.
A coincidence $E_{y_0}(H')=E'_{x_i}$ can happen with
 only one vertex $y_0\in H$ and only one index $i$ ($1\le i\le \ell$),
 namely when $|H'\cap E'_{x_i}|=k-1$.
  If this situation occurs, assume that $E_{y_0}(H')=E'_{x_1}$.
Then, since $\cH$ is $(k-1)$-linear,
 for any $y\in H'\smin\{y_0\}$ and for any $i'\ne i$ with $i'>1$
  we have $|E_y(H') \cap E'_{x_{i'}}| \le k-2$, and so $E_{y}(H')$
  meets $E'_{x_2} \cup \,\cdots\, \cup E'_{x_\ell}$ in at most
  $(\ell-1)(k-2) \le (k-1)(k-2)$ vertices,
 one of which is $y_0$.
Therefore the $k-1$ choices of $y\ne y_0$ cover at most
 $(k^3-3k^2+k)-(k^2-3k+1)+1=k^3-4k^2+4k$
 vertices of $E'_{x_2} \cup \,\cdots\, \cup E'_{x_\ell}$.
Hence, the inequality
 $$
   | \left( E'_{x_2} \cup \cdots \cup E'_{x_\ell} \right)
    \smin \bigcup_{y\in H'} E_y(H') | \ge k+1
 $$
 follows by $|X'|\ge n-2k$ and by the assumed lower
 bound on $|E'_{x_2} \cup \,\cdots\, \cup E'_{x_\ell}|$.
Thus, in this case,
 $H'$ can be completed to a member of
 $\cE^0_{k+1}$ in at least $k+1$ different ways.
The situation is even better if
 we have $|E_y(H') \cap E'_{x_i}| \le k-2$ for all $1\le
i\le \ell$.
 Then summing over all $y\in H'$ and all $1\le i\le \ell$,
we obtain the upper bound
$$
  | (H'^* \smin H') \cap (H^* \smin H) | \le k^2 (k-2) ,
$$
so that there are at least
 $|X'\smin H'^*| \ge n - (k^3 - 2k^2 + 3k) \ge k+1$
 ways to extend $H'$ to a member of $\cE^0_{k+1}$
 whenever $n>k^3$ (and $k\ge 2$).
Consequently,
$$
  \left| \cE^0_{k+1} \right| \ge
    \left| {X'\choose k}\smin \cE_k \right|
$$
 holds, and therefore $|\cE_k| + |\cE^0_{k+1}| \ge {n\choose k} -
O(n^{k-1})$ is valid as $n$ gets large, because $|X'| \ge n-2k$.
 \hfill $\Box$

\section{Geometric upper bound}

As we mentioned in the introduction, $n$ points in $\err^d$ may
 generate as few as $n-1\choose d+1$ affine simplexes, each of which
 has more than $d$ points.
Here we show that the number of affine simplexes can be even smaller.

\bp   \label{lower}
 There is an arrangement of\/ $n$ points in\/ $\err^3$, such that
  the number of affine simplexes determined by them is only
  \begin{itemize}
   \item ${n-1\choose 4} - \frac{(n-2)(n-5)}{2}$ \ if\/ $n$ is even,
   \item ${n-1\choose 4} - \frac{(n-3)(n-5)}{2}$ \ if\/ $n$ is odd;
  \end{itemize}
  that is,\/ $\frac{1}{24}n^4-\frac{5}{12}n^3+O(n^2)$.
\ep

\noindent {\it Proof.}\ \
 First, let\/ $n$ be even.
 Take $n-2$ points $x_1, \dots, x_{n-2}$ on a plane $P\subset \err^3$,
  such that no three of them are collinear, moreover all the
  $n/2-1$ lines $\overline{x_{2i-1}x_{2i}}$ are parallel for
  $i=1,2,\dots,n/2-1$.
 Let $x_{n-1}$ and $x_n$ be two points outside $P$,
  such that the line $\overline{x_{n-1}x_{n}}$ is parallel to
  $\overline{x_1x_2}$ (and hence to the other pairs as well).
 We have the following types of affine simplexes:
 \begin{itemize}
  \item quadruples of points in $P$;
  \item quadruples of the form $\{x_{2i-1},x_{2i},x_{n-1},x_n\}$
    ($i=1,2,\dots,n/2-1$);
  \item quintuples of the form $\{x_a,x_b,x_c,x_{n-1},x_n\}$
    ($1\le a<b<c\le n-2$), where $\{2i-1,2i\}\not\subset\{a,b,c\}$
    for any $i$.
 \end{itemize}
  The number of sets of those three types is $n-2\choose 4$,
   $\frac12(n-2)$, and $\frac16(n-2)(n-4)(n-6)$, respectively.

 If $n$ is odd, we take $(n-3)/2$ pairs of points inside $P$ which
  determine lines $\overline{x_{2i-1}x_{2i}}$ parallel to
   $\overline{x_{n-1}x_{n}}$, plus one point $x_{n-2}$ of $P$
  which is not collinear with any two of $x_1,\dots,x_{n-3}$.
 Then we have $n-2\choose 4$ affine simplexes inside $P$, further
  $\frac12(n-3)$ ones of the form $\{x_{2i-1},x_{2i},x_{n-1},x_n\}$,
   moreover
  $\frac16(n-3)(n-5)(n-7)$ of the form $\{x_a,x_b,x_c,x_{n-1},x_n\}$
    ($1\le a<b<c\le n-3$) where $\{2i-1,2i\}\not\subset\{a,b,c\}$,
   and finally $\frac12(n-3)(n-5)$ of the form
    $\{x_a,x_b,x_{n-2},x_{n-1},x_n\}$
   not containing any pair $\{x_{2i-1},x_{2i}\}$.
  \hfill $\Box$

\section{The \ybl\ inequality}

Recall from Section \ref{s:Sp} that $s(n,k)$ and $s'(n,k)$ are defined as the
 minimum of the sum $\sum_{S\in\cS_k(\cH)} {n\choose |S|}^{-1}$ where $\cH$
 runs over all hypergraphs of order $n$ --- with or without assuming
  $(k-1)$-linearity --- and $\cS_k(\cH)=\cE_k\cup\cE^0_{k+1}$.
Here we study the asymptotic behavior of these two functions, and
 point out a relation to Tur\'an numbers.

\subsection{The limits of $s(n,k)$ and $s'(n,k)$}

Our goal in this subsection is to prove the following result.

\thm   \label{Sper-th}
 For every fixed\/ $k\ge 2$, the limits
  $$
    s_k := \lim_{n\to\infty} s(n,k)
      \qquad \mbox{\rm and} \qquad
    s'_k := \lim_{n\to\infty} s'(n,k)
  $$
  exist and satisfy
  $$
    0 < s'_k \le s_k < 1
  $$
   with strict inequality at both ends.
\ethm

We state three assertions below which together will immediately imply
 the validity of the theorem as the middle inequality holds by definition.

\bl
 For every fixed\/ $k$, the sequences\/ $(s(n,k))_{n=k+1}^\infty$
  and $(s'(n,k))_{n=k+1}^\infty$ are non-decreasing.
\el

 \nin {\it Proof. }
For any hypergraph $\cH=(X,\cE)$ on $n$ vertices, let us introduce the
 notation $m_k:=|\cE_k|$ and $m_{k+1}:=|\cE^0_{k+1}|$.
The inequality
$$
  \frac{m_k}{{n\choose k}} + \frac{m_{k+1}}{{n\choose k+1}} \ge b
$$
 is equivalent to
\begin{equation}   \label{limit-eq}
 m_{k+1} \ge b \cdot {n\choose k+1} - \frac{n-k}{k+1} \, m_k
\end{equation}
 for any $b>0$.
We are going to prove that if the analogue of (\ref{limit-eq}) is valid
 for every hypegraph on $n-1$ vertices, then it is valid for $\cH$ as well.
In what follows, assume that it is valid for $n-1$.

 Let $x\in X$ be any vertex.
We derive the hypergraph $\cH_{-x}$ from $\cH$ by removing $x$ from
 all edges and deleting the edges which become smaller than $k$ after
 cutting out $x$.
Note that $\cH_{-x}$ is $q$-linear whenever so is $\cH$
 (no matter what $q$ we choose),
  although the inverse implication is not valid in general.

Let us denote by $m^x_k$ and $m^x_{k+1}$ the values corresponding to
 $m_k$ and $m_{k+1}$ in $\cH_{-x}$, and by $\cE_{k[-x]}$ and $\cE^0_{k+1[-x]}$
 the corresponding families of sets, respectively.
We then have
 $$
   m^x_{k+1} \ge b \cdot {n-1\choose k+1} - \frac{n-k-1}{k+1} \, m^x_k
      \eqno{(\ref{limit-eq}_x)}
 $$
 by assumption.

A $k$-element set $F$ occurs in $\cE_{k[-x]}$ if and only if it is contained
 in some edge of $\cH_{-x}$; and this happens precisely when $x\notin F$ and
 some $E\in\cE$ contains $F$ as a subset.
Thus, each $F\in\cE_k$ gives rise to a member of $\cE_{k[-x]}$ for exactly
 $n-k$ choices of $x$, and no more sets occur in $\bigcup_{x\in X} \cE_{k[-x]}$.
Similarly, each $F\in\cE^0_{k+1}$ yields a member of $\cE^0_{k+1[-x]}$ for
 exactly $n-k-1$ choices of $x$, and these are all the sets in
  $\bigcup_{x\in X} \cE^0_{k+1[-x]}$.
As a consequence, the equalities
 $$
   \sum_{x\in X} m^x_k = (n-k) \cdot m_k \qquad \mbox{\rm and} \qquad
     \sum_{x\in X} m^x_{k+1} = (n-k-1) \cdot m_{k+1}
 $$
 hold.
Thus, summing up (\ref{limit-eq}$_x$) for all $x\in X$ we obtain
 $$
   (n-k-1) \cdot m_{k+1} \ge bn \cdot {n-1\choose k+1}
    - \frac{n-k-1}{k+1} \, (n-k) \cdot m_k
 $$
  which is equivalent to (\ref{limit-eq}).
This completes the proof.
 \hfill $\Box$

\bl
 For every fixed\/ $k$, we have\/ $s(k+1,k)=s'(k+1,k)=\frac{1}{k+1}$.
\el

\nin {\it Proof. }
 The minimum for both $s(k+1,k)$ and $s'(k+1,k)$ is attained by the
  hypergraph with $k+1$ vertices and with precisely one edge of
  cardinality $k$.
\hfill $\Box$

\bl
 For every fixed\/ $k$, we have\/ $s_k\le 1-\frac{k}{2^k}$.
\el

\nin {\it Proof. }
 To simplify notation, we consider $2n$ vertices instead of $n$.
Consider the hypergraph whose edge set $\cE$ consists of just two disjoint
 sets of cardinality $n$ each.
It is linear, of course. Moreover, we clearly have
 $$
   |\cE_k| = 2 {n\choose k}
     \qquad \mbox{\rm and} \qquad
   |\cE^0_{k+1}| = {2n\choose k+1} - 2 {n\choose k+1} - 2 n {n\choose k}
 $$
 because a $(k+1)$-tuple does \emph{not} belong to $\cE^0_{k+1}$
 if and only if either it is contained in one of the two edges or it
 meets one of the edges in precisely one vertex and the
  other edge in $k$ vertices.
Thus,
 \begin{eqnarray}
  s(2n,k) & \le & \frac{2 {n\choose k}}{{2n\choose k}} + 1 -
    \frac1{{2n\choose k+1}} \left( 2 {n\choose k+1} + 2 n {n\choose k} \right)
 \nonumber \\
  & = & 1 - \frac2{{2n\choose k}} \left( \frac{k+1}{2n-k}
    \left( n {n\choose k} + {n\choose k+1} \right) - {n\choose k}
      \right)
 \nonumber \\
  & = & 1 - \frac{{n\choose k}}{{2n\choose k}} \cdot 2 \cdot
    \left( \frac{n(k+1)}{2n-k} + \frac{n-k}{2n-k} - 1 \right)
 \nonumber \\
  & = & 1 - \frac{{n\choose k}}{{2n\choose k}} \cdot \frac{2nk}{2n-k}
 \nonumber
 \end{eqnarray}
The function in the last line clearly tends to $1-\frac{k}{2^k}$
 as $n\to\infty$, therefore $s_k$ cannot be larger.
 \hfill $\Box$

\subsection{Tur\'an numbers}

For a fixed $k$-uniform hypergraph $\cF$, we use the standard notation
 $\ex(n,\cF)$ for its Tur\'an number; that means the maximum number of
 edges in a $k$-uniform hypergraph of order $n$ which does not contain
 any subhypergraph isomorphic to $\cF$.
Further, let $\cK^{(k)}_{k+1}$ denote the hypergraph with $k+1$ vertices
 and $k+1$ edges of $k$ vertices each (i.e., the complete $k$-uniform
 hypergraph of order $k$).
If $k=2$ then $\cK^{(2)}_3$ is just the triangle $K_3$, the complete
 graph of order $3$.
In this very particular case the equality
 $\ex(n,K_3)=\left\lfloor \frac{n^2}{4} \right\rfloor$
  is well known to hold, but for larger $k$ the determination of
  $\ex(n,\cK^{(k)}_{k+1})$ is a famous open problem
  in extremal hypergraph theory
 (see, e.g., \cite{Sid-srv} for a survey and \cite{FurS} for
  many further references).

\brm
 If\/ $\cH=(X,\cE)$ is a\/ $k$-uniform hypergraph of order\/ $n$ such that
  each\/ $(k+1)$-tuple of vertices contains at least one edge of\/ $\cH$, then\/
  $\cE^0_{k+1}=\es$.
 In particular, taking\/ $\cH$ as the complement of a hypergraph extremal
  for\/ $\ex(n,\cK^{(k)}_{k+1})$, we obtain:
  $$
    s'(n,k) \le 1 - \frac{\ex(n,\cK^{(k)}_{k+1})}{{n\choose k}}.
  $$
 As a consequence,
  $$
    s'_k \le 1 - \lim_{n\to\infty} \frac{\ex(n,\cK^{(k)}_{k+1})}{{n\choose k}}
  $$
  where the limit exists for every fixed\/ $k$, as proved in \cite{KNS}.
 Hence, any lower bound on the Tur\'an density of\/ $\cK^{(k)}_{k+1}$
  implies an upper bound on\/ $s'_k$.
\erm

Note that an analogous implication in the
 opposite direction does not work: upper bounds on $\ex(n,\cK^{(k)}_{k+1})$
 do not imply lower bounds on $s'(n,k)$.
On the other hand, applying the results of Sidorenko \cite{Sid} on
 $\ex(n,\cK^{(k)}_p)$, from the case $p=k+1$ we obtain the following inequality:

\bc   \label{c'k}
 For every\/ $k\ge 3$ we have\/  $s'_k\le \left(1-\frac{1}{k}\right)^{k-1}$.
\ec

It cannot be guaranteed in general that the hypergraphs derived from the
 extremal ones for the Tur\'an problem lead to constructions of
 $(k-1)$-linear hypergraphs, hence they cannot automatically imply
 upper bounds on $s_k$.
But this can be done if $k=2$, and the following exact formula is valid.

\thm
 For every\/ $n\ge 3$ we have\/
  $s(n,2) = s'(n,2) =   1 -
    \left\lfloor \frac{n^2}{4} \right\rfloor\!\left/{n\choose 2}\right.$,
  and therefore\/ $s_2=s'_2=1/2$.
 A hypergraph\/ $\cH=(X,\cE)$ is extremal for\/ $s'(n,2)$ if and only if\/
  $\cE_2$ is the complementary graph of the complete bipartite graph
  $K_{\lfloor n/2 \rfloor, \lceil n/2 \rceil}$; and for\/ $s(n,2)$
  the extremal hypergraph is unique up to isomorphism.
\ethm

 \nin {\it Proof. }
For an upper bound, let $\cH=(X,\cE)$ consist of $n$ vertices and two
 vertex-disjoint edges $E_1,E_2$ with $|E_1| = \left\lfloor n/2 \right\rfloor$
 and $|E_2| = \left\lceil n/2 \right\rceil$.
Then $|\cE_2| = {n\choose 2} - \left\lfloor \frac{n^2}{4} \right\rfloor$
 and $\cE^0_3=\es$.
Since $\cH$ is 1-linear,
 the upper bound follows for both $s(n,2)$ and $s'(n,2)$.
From the argument below, it will also turn out that this is the unique
 1-linear hypergraph attaining equality.

 To prove the lower bound, let $\cH=(X,\cE)$ be any hypergraph.
Note that $\cE_2$ is just a graph; we denote its \emph{complement}
 by $G=(X,E)$, i.e.\ an unordered vertex pair $x_ix_j$ belongs to the
  edge set $E$ of $G$ if and only if $\{x_i,x_j\}\notin\cE_2$.
 Then $\cE^0_3$ is the family of triangles ($K_3$-subgraphs) in $G$.
As long as $G$ is triangle-free, we have
 $|\cE_2|\ge {n\choose 2} - \left\lfloor \frac{n^2}{4} \right\rfloor$
 and the lower bound follows for $\cH$, with equality if and only if
 $G\simeq K_{\lfloor n/2 \rfloor, \lceil n/2 \rceil}$.
Assuming that $|\cE_2|$ is smaller, we have $|E|>n^2\!/4$.

We write the number $|E|$ of edges in $G$ in the form
  $m = \frac{n^2}{4} + \ell$; hence $\ell\ge 1$ is an integer
 if $n$ is even, and $\ell + \frac14 \ge 1$ is an integer
 if $n$ is odd.
It is well known that $G$ has at least
 $\frac{4m^2-mn^2}{3n}$ triangles \cite{MM,NS}.
Thus,
 $$
   |\cE^0_3| \ge \frac{4m}{3n} \left( m - \frac{n^2}{4} \right)
     = \left( \frac{n}{3} + \frac{4\ell}{3n} \right) \ell
       > \frac{n \ell}{3}
 $$
  and consequently
 \begin{eqnarray}
   \frac{|\cE_2|}{{n\choose 2}} + \frac{|\cE^0_3|}{{n\choose 3}}
     & > & \frac{{n\choose 2} - \frac{n^2}{4} - \ell}{{n\choose 2}}
       + \frac{ \frac{n\ell}3}{{n\choose 3}}
 \nonumber \\
     & = & \frac{{n\choose 2} - \frac{n^2}{4}}{{n\choose 2}}
       + \frac{ \frac{n\ell}3 -  \frac{(n-2)\ell}3}{{n\choose 3}}
 \nonumber \\
     & \ge & \frac{{n\choose 2} -
          \left\lfloor \frac{n^2}{4} \right\rfloor}{{n\choose 2}}
       + \frac{2\ell}{3 {n\choose 3}} - \frac1{4{n\choose 2}} .
 \nonumber
 \end{eqnarray}
  This proves the stated inequality for all $\ell \ge (n-2)/8$.
Moreover, if $n$ is even, the theorem follows for all $\ell>0$
 because in that case we need not subtract $1/4$ when moving
 from $\frac{n^2}{4}$ to $\left\lfloor \frac{n^2}{4} \right\rfloor$.

For the rest of the proof, we assume that $n$ is odd and
 $0 < \ell < (n-2)/8$.
Since $\ell$ is relatively small, $G$ must contain some vertex
 $x$ of degree at most $\frac{n-1}2$, for otherwise the number of
 edges would be at least $\frac{n(n+1)}4$, yielding the contradiction
  $\ell \ge \frac n4$.
Let now $\ell'=\ell+\frac14$, that means
 $m=\left\lfloor\frac{n^2}{4}\right\rfloor + \ell'$,
  and consider the graph $G':=G-x$.
This $G'$ has $n'=n-1$ vertices and at least
 $$
   m' := \frac{n^2}{4} + \ell - \frac{n-1}2
     = \frac{(n-1)^2+1}{4} + \ell
       = \frac{(n')^2}{4} + \ell'
 $$
  edges.
Therefore, by the theorem cited above, $G'$ contains at least
 $$
   \frac{4m'}{3n'} \left( m' - \frac{(n')^2}{4} \right)
     = \left( \frac{n-1}{3} + \frac{4\ell'}{3n'} \right) \ell'
       > \frac{(n-1) \ell'}{3}
 $$
  triangles, which certainly is a lower bound on $|\cE^0_3|$, too.
Thus, with a slight modification of the computation above,
 we obtain that
 \begin{eqnarray}
   \frac{|\cE_2|}{{n\choose 2}} + \frac{|\cE^0_3|}{{n\choose 3}}
     & > & \frac{{n\choose 2} - \left\lfloor\frac{n^2}{4}\right\rfloor
      - \ell'}{{n\choose 2}}
       + \frac{ \frac{(n-1)\ell}3'}{{n\choose 3}}
 \nonumber \\
     & = & \frac{{n\choose 2} - \left\lfloor\frac{n^2}{4}\right\rfloor}{{n\choose 2}}
       + \frac{ \frac{(n-1)\ell'}3 - \frac{(n-2)\ell'}3}{{n\choose 3}}
 \nonumber \\
     & = & \frac{{n\choose 2} -
          \left\lfloor \frac{n^2}{4} \right\rfloor}{{n\choose 2}}
       + \frac{\ell'}{3 {n\choose 3}} .
 \nonumber
 \end{eqnarray}
  This completes the proof of the theorem.
\hfill $\Box$

\section{Concluding remarks}

Motivated by a problem arisen in chemistry/stoichiometry, we established
 asymptotically tight extremal results on geometric point sets
 and on finite set systems.
Below we formulate some problems and conjectures that remain open.

\paragraph{Geometry vs.\ hypergraph theory.}

We proved matching asymptotic lower and upper bounds of the form
  ${n\choose d+1}-\Theta(n^d)$ on the minimum number of affine simplexes,
 for every set of $n$ points in $\err^d$ not containing any
 affine simplex of fewer than $d+1$ points.
Our method was to put the problem in a more general context and
 to estimate an extremal function for a class of hypergraphs
  (called $q$-linear, implying the solution for geometric sets
  when $q=d$).
There remains a gap of order $\Theta(n^d)$, however, between the lower
 and upper bounds.

\begin{problem}   \label{pr:geom}
 Given the integers\/ $n$ and\/ $d$, determine the minimum number of
  affine simplexes generated by\/ $n$ points in\/ $\err^d$, no\/ $d$ of which
  lie on a\/ $(d-2)$-dimensional hyperplane.
\end{problem}

\begin{problem}   \label{pr:hg}
 Given the integers\/ $n$ and\/ $k$, determine the minimum value of
  $$|\cE_k|+|\cE^0_{k+1}|$$
  taken over all\/ $(k-1)$-linear hypergraphs\/ $\cH=(X,\cE)$ on\/
  $n$ vertices.
\end{problem}

\begin{problem}    \label{pr:uneq}
 For which values of\/ $k=d+1$ is the minimum for hypergraphs in
  Problem \ref{pr:hg}  equal to that for
  point sets in\/ $\err^d$ in Problem \ref{pr:geom},
  for all\/ $n>n_0(d)$?
\end{problem}

For the case of $d=2$,
 it was proved in \cite{SzL-EJC1998} that the minimum number of
 affine simplexes in $\err^2$ determined by $n\ge 8$ points
 is attained by placing the points on two lines:
 one of the lines contains $n-2$ of the points and the other line
 contains 3 points (and so their intersection point is also selected).
That is, the minimum for Problem \ref{pr:geom} with $d=2$ is
 $${n-2\choose 3} + {n-3\choose 2} + 1 . $$
The construction implies the same upper bound for Problem \ref{pr:hg}
 with $k=3$, attained by the linear (that is, 2-linear) hypergraph with $n$
 vertices and two edges, one of size $n-2$ and the other of size 3.
Moreover, the proof of the matching lower bound in \cite{SzL-EJC1998} gets through
 for linear hypergraphs as well, since it only applies modifications in the
 incidence structure, without any particular geometric assumptions.
Thus, the minimum is the same for Problem \ref{pr:geom} with $d=2$
 and Problem \ref{pr:hg} with $k=3$.

It is not clear, however, whether the answer to Problem \ref{pr:uneq} is
 positive or negative for $d\ge 3$.
We note that the extremal construction of \cite{SzSz-JMC2011}
 cannot be applied for our problem to derive an upper bound on
 affine simplexes in $\err^3$, because in \cite{SzSz-JMC2011}
 the points are arranged in two (equal or nearly equal) collinear sets.
Nevertheless, the following conjecture looks easier than the exact
 determination of minimum.

\bcj
 For every\/ $k\ge 3$ there exists a hypergraph\/ $\cH$ of order\/ $n$ which is
  extremal for Problem \ref{pr:hg} and has\/ $O(n^{k-3})$ edges as\/ $n$ gets large.
\ecj

We note further that the upper bound on $s'_k$ in Corollary~\ref{c'k}
 tends to zero as $k$ gets large, and at present we do not have
 any geometric constructions with the same property for $s_k$.

\bcj
 There exists an integer\/ $k_0$ such that, for every\/ $k\ge k_0$,
  we have\/ $s'_k<s_k$.
\ecj

Perhaps the guess $k_0=3$ is too brave, but we cannot disprove
 even that at present.


\paragraph{Stoichiometry.}

For the original problem originating from \cite{KP-ICE1985} in stoichiometry,
 our Theorem \ref{alg-d} and Corollary \ref{alg-3} imply:
 \begin{itemize}
  \item
 There are at least ${n\choose d+1}- {O}\left(
n^{d}\right) \approx \frac{n^{d+1}}{\left( d+1\right) !}$
minimal reactions among $n$ species if the species are built up
from $d$ kinds of atoms (atomic particles), and if the number of
species (molecules) forming any minimal reaction must be greater than $d$.
 \item
  For the first case previously unsolved, namely $d=3$, the asymptotically
  tight lower bound is ${n\choose 4}-O\left( n^{3}\right) \approx
\frac{n^{4}}{24}$, if reactions with
  three or fewer species are not possible.
   (Especially parallel species, i.e.\ multiple
doses are also excluded.)
 \end{itemize}

\end{document}